\documentclass[12pt,a4paper]{amsart}
\usepackage{amsfonts,amssymb,amscd,amsmath,enumerate,verbatim,calc}
\usepackage{mathrsfs}
\usepackage[all]{xy}
\headsep=1cm \numberwithin{equation}{section}
\newtheorem{thm}{Theorem}[section]
\newtheorem{cor}[thm]{Corollary}
\newtheorem{lem}[thm]{Lemma}

\newtheorem{prop}[thm]{Proposition}
\theoremstyle{definition}

\theoremstyle{remark}

\numberwithin{equation}{section}

\newcommand\Supp{\operatorname{Supp}}
\newcommand\Ass{\operatorname{Ass}}
\newcommand\Assh{\operatorname{Assh}}
\newcommand\mAss{\operatorname{mAss}}

\newcommand\Ann{\operatorname{Ann}}

\newcommand\Soc{\operatorname{Soc}}
\newcommand\Spec{\operatorname{Spec}}
\newcommand\Rad{\operatorname{Rad}}
\newcommand\cd{\operatorname{cd}}
\newcommand\Hom{\operatorname{Hom}}
\newcommand\Ext{\operatorname{Ext}}

\newcommand\height{\operatorname{height}}
\newcommand\Att{\operatorname{Att}}

\newcommand\m{\operatorname{\frak m}}
\newcommand\n{\operatorname{\frak n}}

\newcommand\p{\operatorname{\frak p}}
\newcommand\q{\operatorname{\frak q}}

\begin{document}\title[Cofiniteness over Noetherian complete local rings ]{Cofiniteness over Noetherian complete local rings }
\author[K. Bahmanpour ]{  Kamal Bahmanpour  }

\address{Department of Mathematics, Faculty of Sciences, University of Mohaghegh Ardabili,
56199-11367, Ardabil, Iran;
and School of Mathematics, Institute for Research in Fundamental Sciences (IPM),
 P.O. Box. 19395-5746, Tehran, Iran.} \email{\it bahmanpour.k@gmail.com}

\thanks{ 
This research of the author was supported by a grant from IPM (No. 96130018).}

\begin{abstract}
   In this paper we prove the following generalization of a result of
Hartshorne: Let $(S,\n)$ be a regular local ring of dimension $4$. Assume that $x,y,u,v$ is a regular system of parameters for $S$ and $a:=xu+yv$. Then for each finitely generated $S$-module $N$ with $\Supp N=V(aS)$ the socle of  $H^2_{(u,v)S}(N)$ is infinite dimensional. Also, using this result, for any commutative Noetherian complete local  ring $(R,\m)$, we characterize the class of all ideals $I$ of $R$ with the property that, for every finitely generated $R$-module $M$, the local cohomology modules $H^i_I(M)$ are $I$-cofinite for all $i\geq 0$.
  \end{abstract}
\maketitle
\section{Introduction}

Throughout this paper, let $R$ denote a commutative Noetherian ring
(with identity) and $I$ be an ideal of $R$. For an $R$-module $M$, the
$i$th local cohomology module of $M$ with support in $V(I)$
is defined as:
$$H^i_I(M) = \underset{n\geq1} {\varinjlim}\,\, \Ext^i_R(R/I^n,
M).$$   We refer the reader to \cite{BS} or \cite{Gr1} for more
details about local cohomology.\\

Recall that, if $(R,\m,k)$ is a local ring then for any $R$-module $L$
{\it the socle of $L$} denoted by $\Soc_{R} L$ is defined as $ \Soc_{R} L=(0:_L\m)\cong \Hom_R(k ,L) $, which is a $k$-vector space.\\

In Ref. \cite{Hu} Huneke conjectured the following:\\

{\bf Conjecture:} {\it For any ideal $I$ in a regular local ring
$(R, \frak m)$, the $R$-module $\Soc_R H^i_I(R)$ is finitely generated for each integer $i\geq 0$.} \\

 It is shown by Huneke and Sharp \cite{HS} and Lyubeznik \cite{Ly1,Ly2} that this
conjecture holds for any regular local ring containing a field. The first example of a
 local cohomology module with an infinite dimensional socle was constructed by Hartshorne. \\

{\bf Hartshorne's example:} (See \cite[\S3]{Ha}) {\it Let $k$ be a field, $R = k[[u, v]][x, y]$, $I = (x, y)R$,
$P = (u, v, x, y)R$ and $f = ux + vy$. Then $\Soc_{R_P} H_{IR_P}^2(R_{P}/fR_{P})$ is
infinite dimensional.}\\

 Hartshorne proved this by exhibiting an infinite set of linearly independent
elements in the $\Soc_{R_P}  H_{IR_P}^2(R_{P}/fR_{P})$. In 2004 a similar family of
such informative examples was constructed by Marley and Vassilev in Ref. \cite{MV1}. Beyond that work, however, there are a few results in the literature which explain or
generalize Hartshorne's example.\\

 In section 2 of this paper, as a  generalization of the Harthshorne's example, we shall prove the following theorem:\\

 {\bf  Theorem 1.} {\it  Let $(S,\n)$ be a regular local ring of dimension $4$. Assume that $x,y,u,v$ is a regular system of parameters for $S$ and $a:=xu+yv$. Then for each finitely generated $S$-module $N$ with $\Supp N=V(aS)$, the $S$-module  $H^2_{(u,v)S}(N)$ is of dimension zero with infinite dimensional socle. }\\

For an $R$-module $M$, the notion $\cd (I, M)$, {\it the
 cohomological dimension of $M$ with respect to $I$}, is defined as:
$$\cd(I,M)={\rm sup}\{i\in \Bbb{N}_0\,\,:\,\,H^i_I(M)\neq 0\}$$
and the notion $q(I,M)$, which for first time was introduced by Hartshorne, is defined as:
 $$q(I,M)={\rm sup}\{i\in \Bbb{N}_0\,\,:\,\,H^i_I(M)\,\,{\rm is }\,\,{\rm not }\,\,{\rm Artinian } \},$$  with the usual convention that the
supremum of the empty set of integers is interpreted as $-\infty$. These two notions have been studied by
several authors, (see \cite{B, DNT, DY, Fa, GBA, Ha2, HL}). \\

Hartshorne in Ref. \cite{Ha}  defined an $R$-module $L$ to be
$I$-{\it cofinite}, if $\Supp L\subseteq
V(I)$ and ${\rm Ext}^{i}_{R}(R/I, L)$ is a finitely generated module
for all $i$. Then he posed the following question:\\

{\bf Question 1:} {\it For which Noetherian rings $R$ and ideals $I$ are the modules
$H^{i}_{I}(M)$ $I$-cofinite for  all finitely generated
$R$-modules $M$ and all $i\geq0$}?\\

In this paper, we denote by $\mathscr{I}(R)$ the class of all ideals $I$ of $R$
with the property that, for every finitely generated $R$-module $M$, the local cohomology modules $H^i_I(M)$ are $I$-cofinite for all $i\geq 0$. \\

Concerning the Question 1, there are several papers in the literature containing some sufficient conditions for the ideals of $R$ being in $\mathscr{I}(R)$, (see \cite{B2, B, BN, Ch, De, DM, HK, Ka0, MV, Yo}). In section 3 of this paper, in order to finding a necessary and sufficient condition for the ideals of any Noetherian complete local ring $(R,\m)$ being in $\mathscr{I}(R)$, we will focus on the results of the more recently published article \cite{B2}. \\

We recall that the present author in Ref. \cite{B2}, for any ideal $I$ of $R$ and any finitely generated $R$-module $M$, defined: $$\mathfrak{A}(I,M):=\{\p \in \mAss_R M\,\,:\,\,I+\p=R\,\,\,\,{\rm or}\,\,\,\,\p\supseteq I\},$$$$\mathfrak{B}(I,M):=\{\p \in \mAss_R M\,\,:\,\,\cd(I,R/\p)=1\},$$$$\mathfrak{C}(I,M):=\{\p \in \mAss_R M\,\,:\,\,q(I,R/\p)=1\}\,\,{\rm and}$$$$\mathfrak{D}(I,M):=\{\p \in \mAss_R M\,\,:\,\,0\leq \dim R/(I+\p)\leq 1\}.$$

Also, he proved that if $\mAss_R R=\mathfrak{A}(I,R)\cup\mathfrak{B}(I,R)\cup \mathfrak{C}(I,R)\cup\mathfrak{D}(I,R)$ then $I\in \mathscr{I}(R)$.  Then, he posed the following two questions:\\

{\bf Question A:} {\it Let $(R,\m)$ be a Noetherian complete local ring and $I\in\mathscr{I}(R)$. Whether $\mAss_R R= \mathfrak{A}(I,R) \cup\mathfrak{B}(I,R) \cup\mathfrak{D}(I,R)$?} \\

{\bf Question B:} {\it Let $(R,\m)$ be a Noetherian complete local domain of dimension $4$ and $I$ be an ideal of $R$ with $\height I=2$. Whether $I\not\in\mathscr{I}(R)$?} \\

 He proved that Question {\rm A} has an affirmative answer in general if and only if Question {\rm B} has so, (see \cite[Proposition 4.7]{B2}).\\

In section 3 of this paper, using the results of section 2, we present an affirmative answer to Question {\rm B}. Then, we specify the elements of $\mathscr{I}(R)$, when $(R,\m)$ is a Noetherian complete local ring. More precisely, we prove the following theorem:\\

{\bf Theorem 2.} {\it Let $(R,\m)$ be a Noetherian complete local ring and $I$ be an ideal of $R$. Then the following statements are equivalent:}
\begin{enumerate}[\upshape (i)]
  \item  $I\in\mathscr{I}(R)$.
  \item  $\mAss_R R= \mathfrak{A}(I,R) \cup\mathfrak{B}(I,R) \cup\mathfrak{D}(I,R)$.
  \end{enumerate}

Also, as an immediate consequence of Theorem 2, we deduce the following result:\\

{\bf  Corollary. } {\it Let $(S,\n)$ be a regular complete local ring. Then}
  $$\mathscr{I}(S)=\{I\leq S\,\,:\,\,\Rad(I)=xS,\,\,{\rm for}\,\,{\rm some}\,\,x\in S\}\bigcup\{I\leq S\,\,:\,\, 0\leq\dim S/I\leq 1\}.$$

Throughout this paper, for each $R$-module $L$, we denote by
 $\Assh_R L$ (respectively by $\mAss_R L$), the set $\{\p\in \Ass_R L\,\,:\,\, \dim R/\p= \dim L\}$  (respectively the set of
 minimal elements of $\Ass_R L$ with respect to inclusion). For any ideal $\frak{a}$ of $R$, we denote
$\{\frak p \in {\rm Spec}\,R:\, \frak p\supseteq \frak{a}\}$ by
$V(\frak{a})$. For any ideal $\frak{b}$ of $R$, {\it the
radical of} $\frak{b}$, denoted by $\Rad(\frak{b})$, is defined to
be the set $\{x\in R \,: \, x^n \in \frak{b}$ for some $n \in
\mathbb{N}\}$.  For any unexplained notation and terminology we refer the reader to \cite{BS} and \cite{Mat}.

\section{A new lookout to Hartshorne's example}

The main purpose of this section is to prove Theorem 2.7, which is a generalization of an incredibly rare and valuable example constructed by Hartshorne in \cite[\S3]{Ha}.\\

 We start this section with some auxiliary lemmas, which are needed in the proofs of Proposition 2.6 and Theorem 2.7.

 \begin{lem}
 \label{2.1}
Let $(S,\n,k)$ be a regular local ring of dimension $4$ and $x,y,u,v$ be a regular system of parameters for $S$. Then $a:=xu+yv$ is a prime element of $S$.
 \end{lem}
\proof  From the hypothesis we know that $S$ is a unique factorization domain and $a$ is neither 0 nor a unit; so, in order to prove the assertion, it is enough to prove that $a$ is an irreducible element of $S$. In contrary, assume that $a$ is not irreducible.  Then $a$ can be expressed as a product of two elements of $\n$. Suppose that $a=bc$ with $b,c\in \n$. Then there are elements $\alpha_i, \beta_i\in S$, for $1\leq i \leq 4$, such that $$b=\alpha_1x+\alpha_2y+\alpha_3u+\alpha_4v\,\,\,\,{\rm and}\,\,\,\,c=\beta_1x+\beta_2y+\beta_3u+\beta_4v.$$Also, as $a\in \n^2$ and $a\not\in \n^3$ we have $b\not\in \n^2$ and $c\not\in\n^2$. Therefore, each of the sets $\{\alpha_1,\alpha_2,\alpha_3,\alpha_4\}$ and $\{\beta_1,\beta_2,\beta_3,\beta_4\}$ contains at least one unit. Using the symmetry of the problem, without loss of generality, we may assume that $\alpha_1$ is a unit. Since, $x,y,u,v$ is a regular system of parameters for $S$ it follows that the set $$B:=\{x^2+\n^3,y^2+\n^3,u^2+\n^3,v^2+\n^3,xy+\n^3,xu+\n^3,xv+\n^3,yu+\n^3,yv+\n^3,uv+\n^3\}$$
is a base for the $k$-vector space $\n^2/\n^3$. From the relation $$xu+yv+\n^3=a+\n^3=bc+\n^3=(\alpha_1x+\alpha_2y+\alpha_3u+\alpha_4v)(\beta_1x+\beta_2y+\beta_3u+\beta_4v)+\n^3 $$
we get the relation
$$\big[(\alpha_1\beta_1)x^2+(\alpha_2\beta_2)y^2+(\alpha_3\beta_3)u^2+(\alpha_4\beta_4)v^2+
(\alpha_1\beta_2+\alpha_2\beta_1)xy+(\alpha_1\beta_3+\alpha_3\beta_1-1_{_S})xu+$$
$$(\alpha_1\beta_4+\alpha_4\beta_1)xv+(\alpha_2\beta_3+\alpha_3\beta_2)yu+(
\alpha_2\beta_4+\alpha_4\beta_2-1_{_S})yv+(\alpha_3\beta_4+\alpha_4\beta_3)uv\big]+\n^3=\n^3.$$
Since, $B$ is a set of linearly independent vectors over the field $k$, we get the following relations:\\

  \item (1) $\alpha_1\beta_1\in \n$,
  \item (2) $\alpha_1\beta_2+\alpha_2\beta_1\in\n$,
  \item (3) $\alpha_1\beta_4+\alpha_4\beta_1\in\n$,
  \item (4) $\alpha_2\beta_4+\alpha_4\beta_2-1_{_S}\in\n$.\\

As $\alpha_1$ is a unit, by the relation (1) we have $\beta_1\in \n$. Therefore, the relations (2) and (3) imply that $\beta_2, \beta_4\in \n$; so, the relation (4) yields that $1_{_S}\in \n$, which is a contradiction.\qed\\

\begin{lem}
\label{2.2}
${\rm (See}$ \cite[Theorem 4.9]{B}$)$
Let $R$ be a Noetherian ring and $I$ be an ideal of $R$. Then for each finitely generated $R$-module $M$ with $q(I,M)\leq1$, the $R$-modules $H^i_I(M)$ are $I$-cofinite for all $i\geq0$.
\end{lem}\qed\\

\begin{lem}
\label{2.3}
${\rm (See}$ \cite[Lemma 2.3]{B0}$)$
Let $(R,\m)$ be a Noetherian complete local ring and $I$ be an ideal of $R$. If ${\rm cd}(I,R)=t$ and the $R$-module $H^t_I(R)$ is Artinian and $I$-cofinite, then $$\Att_R H^t_I(R)=\{\p \in \mAss_R R\,\,:\,\,\dim R/\p =t\,\,{\rm and}\,\,\Rad(I+\p)=\m\}.$$
\end{lem}\qed\\

\begin{lem}
\label{2.4}
${\rm (See}$ \cite[Theorem 2.2]{DNT}$)$
Let $R$ be a Noetherian ring, $I$ be an ideal of $R$ and $M$, $N$ be two finitely generated $R$-modules. If
$\Supp M \subseteq \Supp N$ then $\cd(I,M)\leq \cd(I,N).$
\end{lem}\qed\\

\begin{lem}
\label{2.5}
${\rm (See}$ \cite[Theorem 3.2]{DY}$)$
Let $R$ be a Noetherian ring, $I$ be an ideal of $R$ and $M$, $N$ be two finitely generated $R$-modules. If
$\Supp M \subseteq \Supp N$ then $q(I,M)\leq q(I,N).$
\end{lem}\qed\\

The following proposition plays a key role in the proof of Theorem 2.7.

\begin{prop}
 \label{2.6}
  Let $(S,\n,k)$ be a regular local ring of dimension $4$. Assume that $x,y,u,v$ is a regular system of parameters for $S$ and $a:=xu+yv$. Then the following statements hold:
 \begin{enumerate}[\upshape (i)]
  \item  $\cd((u,v)S,S/aS)= 2$.
  \item  $\Supp H^2_{(u,v)S}(S/aS)=\{\n\}$.
  \item  $\dim_{k} \Soc_S H^2_{(u,v)S}(S/aS)=\infty$.
  \end{enumerate}
\end{prop}

\proof (i) As the ideal $(u,v)S$ of $S$ is generated by two elements, it follows from \cite[Theorem 3.3.1]{BS} that $\cd((u,v)S,S/aS)\leq 2$. Set $M:=S/(x,y)S$. Then, $M$ is a finitely generated $S$-module of dimension $2$. Therefore, using  the {\it Grothendieck's Non-vanishing Theorem} and \cite[Exercise 2.1.9]{BS} we have
$$H^2_{(u,v)S}(M) \cong H^2_{(u,v)S+\Ann_S M}(M)=H^2_{\n}(M)\neq0.$$ Furthermore, as $a\in (x,y)S$ it follows that $\Supp M\subseteq \Supp S/aS$ and so by Lemma 2.4 we have $$\cd((u,v)S,S/aS)\geq \cd((u,v)S,M)\geq 2.$$
Hence, $\cd((u,v)S,S/aS)= 2$.\\

(ii) By part (i) we have $H^2_{(u,v)S}(S/aS)\neq0$ and so $\Supp H^2_{(u,v)S}(S/aS)\neq\emptyset$. Therefore, in order to prove (ii), it is enough to prove $\Supp H^2_{(u,v)S}(S/aS)\subseteq\{\n\}$. In contrary, assume that $\Supp H^2_{(u,v)S}(S/aS)\not\subseteq\{\n\}$. Then, there is an element $\p\in \Supp H^2_{(u,v)S}(S/aS)$ such that $\p\neq \n$. It is clear that $\dim S_{\p}=\height \p<\height \n =\dim S=4$ and hence $\height \p\leq 3$. Moreover, as $$H^2_{(u,v)S_{\p}}(S_{\p}/aS_{\p})\cong (H^2_{(u,v)S}(S/aS))_{\p}\neq0,$$the {\it Grothendieck's Vanishing Theorem} yields that $\dim S_{\p}/aS_{\p}\geq 2$. Thus, using the fact that $S$ is a domain and $a\neq0$, we get $$\height \p=\dim S_{\p}=1+\dim S_{\p}/aS_{\p}\geq 3.$$ So, $\height \p=3$. On the other hand, we have $$\p \in \Supp H^2_{(u,v)S}(S/aS) \subseteq V((u,v)S)$$ and so $(u,v)S\subseteq \p$. Since, by our standard hypothesis, $x,y,u,v$ is a regular system of parameters for $S$, we deduce that the local ring $S/(u,v)S$ is regular too. Hence, $S/(u,v)S$ is a unique factorization domain. Therefore, the relation $$\height \p/(u,v)S=\height \p - \height (u,v)S=3-2=1$$ implies that the ideal $\p/(u,v)S$ of $S/(u,v)S$ is principal. Thus, there exists an element $z\in \p$ such that $\p=(u,v,z)S$. Since, $$(u,v)S\subseteq \p \subset \n=(x,y,u,v)S$$ we can see that at least one of the relations $x\not\in \p$ and $y\not\in \p$ holds.
 Using the symmetry of the problem, without loss of generality, we may assume that $x\not\in \p$. Then, $x/1_{_S}\in S_{\p}$ is a unit and so $(a,v)S_{\p}=(u,v)S_{\p}$. Consequently, we have $$(a,v,z)S_{\p}=(u,v,z)S_{\p}=\p S_{\p}$$ and so $a/1_{_S},v/1_{_S},z/1_{_S}$ is a regular system of parameters for the regular local ring $S_{\p}$. Now, it is clear that $S_{\p}/aS_{\p}$ is a regular local ring of dimension $2$. Hence, the relation $$H^2_{(u,v)S_{\p}/aS_{\p}}(S_{\p}/aS_{\p})\cong H^2_{(u,v)S_{\p}}(S_{\p}/aS_{\p})\neq0,$$ considering the {\it Lichtenbaum-Hartshorne Vanishing Theorem} yields that $$\Rad((u,v)S_{\p}/aS_{\p})=\p S_{\p}/aS_{\p}.$$ Therefore, $(u,v)S_{\p}=\Rad((u,v)S_{\p})=\p S_{\p}$, which is a contradiction.\\

 (iii) In contrary, assume that $\dim_{k} \Soc_S H^2_{(u,v)S}(S/aS)<\infty$. Then, as by part (ii) we have $\Supp H^2_{(u,v)S}(S/aS)=\{\n\}$, it is clear that the $S$-module $H^2_{(u,v)S}(S/aS)$ is Artinian. Consequently, using part (i) we have $q((u,v)S, S/aS)\leq 1$. Therefore, in view of Lemma 2.2, the $S$-module $H^2_{(u,v)S}(S/aS)$ is $(u,v)S$-cofinite. So, the $\widehat{S}$-module $$H^2_{(u,v)S}(S/aS)\otimes_S \widehat{S}\cong H^2_{(u,v)\widehat{S}}(\widehat{S}/a\widehat{S})$$ is Artinian and $(u,v)\widehat{S}$-cofinite, where $\widehat{S}$ is the $\n$-adic completion of $S$. But, as $\widehat{S}$ is a regular local ring of dimension $4$ with the maximal ideal $\n\widehat{S}=(x,y,u,v)\widehat{S}$ it is clear that $x,y,u,v$ is a regular system of parameters for $\widehat{S}$. Thus, by Lemma 2.1, $a$ is a prime element of $\widehat{S}$ and hence $a \widehat{S}$ is a prime ideal of $\widehat{S}$.  Since, $H^2_{(u,v)\widehat{S}}(\widehat{S}/a\widehat{S})$ is Artinian and $(u,v)\widehat{S}$-cofinite it follows that $$H^2_{(u,v)\widehat{S}/a\widehat{S}}(\widehat{S}/a\widehat{S}) \cong H^2_{(u,v)\widehat{S}}(\widehat{S}/a\widehat{S})$$ is an Artinian $(u,v)\widehat{S}/a\widehat{S}$-cofinite $\widehat{S}/a\widehat{S}$-module. Furthermore, by part (i) and the {\it Independence Theorem} we have  $$\cd((u,v)\widehat{S}/a\widehat{S},\widehat{S}/a\widehat{S})= \cd((u,v)\widehat{S},\widehat{S}/a\widehat{S})=2.$$ Now, by Lemma 2.3 we can deduce that $$\Rad((u,v)\widehat{S}/a\widehat{S})=\n \widehat{S}/a \widehat{S},$$ which is a contradiction. \qed\\

Now, we are ready to state and prove the main result of this section.

\begin{thm}
 \label{2.7}
Let $(S,\n,k)$ be a regular local ring of dimension $4$. Assume that $x,y,u,v$ is a regular system of parameters for $S$ and $a:=xu+yv$. Then for each finitely generated $S$-module $N$ with $\Supp N=V(aS)$, the following statements hold:
 \begin{enumerate}[\upshape (i)]
  \item  $\cd((u,v)S,N)= 2$.
  \item  $\Supp H^2_{(u,v)S}(N)=\{\n\}$.
  \item  $\dim_{k} \Soc_S H^2_{(u,v)S}(N)=\infty$.
  \end{enumerate}
\end{thm}
\proof (i) Follows from Proposition 2.6(i) and Lemma 2.4.\\

(ii) From the part (i) it follows that $\Supp H^2_{(u,v)S}(N)\neq \emptyset$. Also, using Lemma 2.4 and localization it follows from parts (i) and (ii) of Proposition 2.6 that  $\Supp H^2_{(u,v)S}(N)\subseteq\{\n\}$. \\

(iii) In contrary, assume that  $\dim_{k} \Soc_S H^2_{(u,v)S}(N)<\infty$. Then, it follows from part (ii) that the $S$-module $H^2_{(u,v)S}(N)$ is Artinian. But, in this situation part (i) implies that $q((u,v)S,N)\leq 1$. Therefore, using the relation $\Supp N=V(aS)=\Supp S/aS$ it follows from Lemma 2.5 that $q((u,v)S,S/aS)\leq 1$. So, the $R$-modules $H^2_{(u,v)S}(S/aS)$ is Artinian and hence $\dim_{k} \Soc_S H^2_{(u,v)S}(S/aS)<\infty$. But by part (iii) of Proposition 2.6, this is a contradiction.\qed\\

\section{A characterization of $\mathscr{I}(R)$ over complete local rings}

The main purpose of this section is to prove Proposition 3.4 and Theorem 3.8. In fact, Proposition 3.4 presents an affirmative answer to a question
raised by the present author in Ref. \cite{B2} and the proof of Theorem 3.8 relies heavily on this result. The following auxiliary lemmas are quite useful in the proof of Proposition 3.4.

\begin{lem}
\label{3.1}
${\rm (See}$ \cite[Corollary 4.4]{B2}$)$
 Let $(R,\m)$ be a Noetherian complete local ring  and $I\in\mathscr{I}(R)$.  Then $\cd(I,R/\p)=\height (I+\p)/\p$, for each $\p \in \Spec R$.
\end{lem}
\qed\\

\begin{lem}
\label{3.2}
${\rm (See}$ \cite[Proposition 4.6]{B2}$)$
 Let $(R,\m)$ be a Noetherian complete local ring  and $I\in\mathscr{I}(R)$.  Then, $\mAss_R R/(\p+I)=\Assh_R R/(\p+I)$, for each $\p\in \Spec R$.
 \end{lem}
\qed\\

\begin{lem}
\label{3.3}
${\rm (See}$ \cite[Corollary 2.7]{BN}$)$
 Let $R$ be a Noetherian ring and $I$ be an ideal of $R$. If $M$ is a finitely generated $R$-module
such that $\dim M/IM\leq1$ then the
$R$-modules $H^{i}_{I}(M)$ are $I$-cofinite for all $i\geq0$.
 \end{lem}
\qed\\

The following proposition is the first main result of this section.

\begin{prop}
\label{3.4}
Let $(R,\m,k)$ be a Noetherian complete local domain of dimension $4$ and $I$ be an ideal of $R$ with $\height I=2$. Then, $I\not\in\mathscr{I}(R)$.
 \end{prop}
\proof If $\cd(I,R)\neq \height I$ then the assertion holds by Lemma 3.1. So, we may assume that $\cd(I,R)=\height I=2$. Also, if $\mAss_R R/I\neq \Assh_R R/I$ then the assertion follows from Lemma 3.2. Therefore, we may assume that $\height \p=2$, for each $\p\in \mAss_R R/I$. Now, in contrary assume that $I\in\mathscr{I}(R)$. Note that as $R$ is a catenary domain we have $\dim R/I = \dim R - \height I=4-2=2$.\\

At the first step of the proof, we deal with the structure theory of complete local rings. But, in order to develop our strategy, first we need to consider each of the following three possible cases separately. In fact, in each of these cases we choose a suitable  element $\q_1\in \mAss_R R/I$ and an appropriate system of parameters $x,y,u,v$ for $R$ with the following other two extra properties:\\

\item (i) $(x,y)R\subseteq \q_1$,
\item (ii) $u,v$ is a system of parameters for the $R$-module $R/I$.\\

Case 1. Assume that $(R,\m,k)$ is equicharacteristic. Pick elements $x,y\in I$ and $u,v\in\m$ such that $x,y,u,v$ is a system of parameters for $R$. Also, pick an arbitrary element $\q_1\in \mAss_R R/I$.  \\

Case 2. Assume that $R$ is of characteristic 0, ${\rm char}\,k=p$ ($p$ a prime integer), and $$p1_{_R}\in \bigcup_{\q\in \mAss_R R/I}\q.$$ Pick  $\q_1\in \mAss_R R/I$ with the property $p1_{_R}\in \q_1$, and set $x:=p1_{_R}$. Then, $\height \q_1=2$ and by the {\it Principal Ideal Theorem} we have $\height \p=1$, for each $\p \in \mAss_R R/xR$. Hence, $\q_1\not\subseteq \bigcup_{\p\in \mAss_R R/xR}\p$ and so we can find an element $y\in \q_1$ with $y\not\in \bigcup_{\p\in \Assh_R R/xR}\p$.\\

 Next, pick an element $u\in\m$ with $$u\not\in\left( \left(\bigcup_{\q\in \Assh_R R/I}\q\right)\bigcup\left(\bigcup_{\p\in \Assh_R R/(x,y)R}\p\right)\right).$$ Finally, we can find an element $v\in \m$ such that $$v\not\in\left( \left(\bigcup_{\q\in \Assh_R R/(I+uR)}\q\right)\bigcup\left(\bigcup_{\p\in \Assh_R R/(x,y,u)R}\p\right)\right).$$ \\

Case 3. Assume that $R$ is of characteristic 0, ${\rm char}\,k=p$ ($p$ a prime integer), and $$p1_{_R}\not\in \bigcup_{\q\in \mAss_R R/I}\q.$$ Set $u:=p1_{_R}$ and pick $v\in \m$ with the property $$v\not \in \left( \left(\bigcup_{\p\in \Assh_R R/uR}\p\right)\bigcup\left(\bigcup_{\q\in \Assh_R R/(I+uR)}\q\right)\right).$$ Then, it is clear that $u,v$ is a system of parameters for the $R$-module $R/I$ and $$\dim R/(u,v)R=2.$$ Thus, we have $\Rad((I+(u,v)R)/I)=\m/I$ and so $\Rad(I+(u,v)R)=\m$. Therefore, there are elements $x,y\in I$ and $z,t\in (u,v)R$ such that $x_1=x+z$, $y_1=y+t$ is a system of parameters for the $R$-module $R/(u,v)R$. Then, we have $\Rad((x_1,y_1)R+(u,v)R)=\m$ and $(x_1,y_1,u,v)R=(x,y,u,v)R$, which means $x,y,u,v$ is a system of parameters for $R$. Now, pick an arbitrary element $\q_1\in \mAss_R R/I$. \\

 By \cite[Theorem 28.3]{Mat}, in the case 1, $R$ contains a coefficient field $F$. Also, in view of \cite[Theorem 29.3]{Mat} in both of the cases 2 and 3, $R$ contains a coefficient ring $A$, such that $A$ is a complete DVR with the maximal ideal $pA$. In the case 1, set $S:=F[[x,y,u,v]]$, in the case 2, set $S:=A[[y,u,v]]$ and in the case 3, set $S:=A[[x,y,v]]$. Then by the proof of \cite[Theorem 24.9]{Mat}, $(S,\n)$ is a complete regular local ring of dimension $4$ and $R$ is finitely generated as an $S$-module, where $\n=(x,y,u,v)S$. In particular, $S\subseteq R$ is an integral extension and $x,y,u,v$ is a regular system of parameters for $S$. \\

Next, set $a:=xu+yv$. Then, by Lemma 2.1, $a$ is a prime element of $S$ and so $aS$ is a prime ideal of $S$. Since, $$aR\subseteq (x,y)R\subseteq \q_1$$ it follows that $\q_1$ contains an element $\q_2\in \mAss_R R/aR$. By the {\it Principal Ideal Theorem} we have $\height \q_2=1$. Since, $S\subseteq R$ is an integral extension and $R$ is a catenary domain it follows that $$\dim S/aS=3=\dim R/\q_2=\dim S/(S\cap \q_2).$$ Furthermore, $aS\subseteq aR\cap S \subseteq \q_2 \cap S$ and both of the ideals $aS$ and $\q_2 \cap S$ are prime. Now, we are ready to deduce that $\q_2 \cap S=aS$. Therefore, $R/\q_2$ is a finitely generated $S$-module with $\Ann_S R/\q_2 = aS$ and so $\Supp S/aS= \Supp R/\q_2$. Now, it follows from Theorem 2.7 that $\Supp H^2_{(u,v)S}(R/\q_2)=\{\n\}$. Therefore, $\Ass_S H^2_{(u,v)S}(R/\q_2)=\{\n\}$. Moreover, by the {\it Independence Theorem} we have $$H^2_{(u,v)R}(R/\q_2)\cong H^2_{(u,v)S}(R/\q_2)\neq 0.$$ In particular, $\Ass_R H^2_{(u,v)R}(R/\q_2)\neq \emptyset$. We claim that $\Ass_R H^2_{(u,v)R}(R/\q_2)=\{\m\}$. To prove this assertion, it is sufficient for us to show that $\Ass_R H^2_{(u,v)R}(R/\q_2)\subseteq\{\m\}$. Let $$\p \in \Ass_R H^2_{(u,v)R}(R/\q_2).$$ Then, $\p\in \Spec R$ and there is an element $h\in  H^2_{(u,v)R}(R/\q_2)$ such that $(0:_Rh)=\p$. Hence, $$(0:_Sh)=(0:_Rh)\cap S=\p\cap S\in \Spec S$$ and so $$\p\cap S\in \Ass_S H^2_{(u,v)R}(R/\q_2)=\Ass_S H^2_{(u,v)S}(R/\q_2)=\{\n\}.$$ Thus, $\p\cap S=\n$. But, as $R$ is integral over $S$, it follows that $\m$ is the one and only prime ideal of $R$ which has contraction to $S$ equal to $\n$. Therefore, $\p=\m$ and hence $\Ass_R H^2_{(u,v)R}(R/\q_2)=\{\m\}$, which implies that $\Supp_R H^2_{(u,v)R}(R/\q_2)=\{\m\}$.\\

At this point, we prove that $\dim_{k} \Soc_R H^2_{(u,v)R}(R/\q_2)=\infty$.  Let us, in contrary, assume that $$\dim_{k} \Soc_R H^2_{(u,v)R}(R/\q_2)<\infty.$$ Then, as $\Supp_R H^2_{(u,v)R}(R/\q_2)=\{\m\}$ it follows that the $R$-module $H^2_{(u,v)R}(R/\q_2)$ is Artinian. In particular, using   \cite[Theorem 3.3.1]{BS} we have $q((u,v)R, R/\q_2)\leq 1$. Therefore, in view of Lemma 2.2, the $R$-module $H^2_{(u,v)R}(R/\q_2)$ is $(u,v)R$-cofinite. So, $$H^2_{((u,v)R+\q_2)/\q_2}(R/\q_2)\cong H^2_{(u,v)R}(R/\q_2)$$ is $((u,v)R+\q_2)/\q_2$-cofinite. Moreover, by \cite[Theorem 3.3.1]{BS} and the {\it Independence Theorem} we have  $$\cd(((u,v)R+\q_2)/\q_2,R/\q_2)= 2.$$ Now, by Lemma 2.3 we can deduce that $$\Rad(((u,v)R+\q_2)/\q_2)=\m/\q_2.$$ But, by the {\it Generalized Principal Ideal Theorem} this implies that $3=\height \m/\q_2\leq 2$, which is a contradiction.\\

Now, we claim that $\q_2 \subseteq \bigcup_{\p\in \mAss_R R/(u,v)R}\p$. In contrary, assume that $$\q_2 \not \subseteq \bigcup_{\p\in \mAss_R R/(u,v)R}\p.$$ Then, we have $\dim R/((u,v)R+\q_2)<\dim R/(u,v)R=2$ and so $\dim R/((u,v)R+\q_2)\leq 1$. Hence, in view of Lemma 3.3, the $R$-module $H^2_{(u,v)R}(R/\q_2)$ is $(u,v)R$-cofinite, which is a contradiction. So, there exists $\p_1\in \mAss_R R/(u,v)R$ such that $\q_2\subseteq \p_1$.\\

Note that as $R$ is a catenary domain we have $\height \p_1 + \dim R/\p_1 =\dim R =4$ and from the fact that $\p_1\in \mAss_R R/(u,v)R$, by the {\it Generalized Principal Ideal Theorem} we have $\height \p_1\leq 2$. Also, it is clear that $\dim R/\p_1 \leq \dim R/(u,v)R=2$. Hence, $\height \p_1 =2 = \dim R/\p_1$. \\

Now, as $u,v$ is a system of parameters for the $R$-module $R/I$ and $(u,v)\subseteq \p_1$ we have$$\m=\Rad(I+(u,v))\subseteq \Rad(I+\p_1)\subseteq \m$$ and so $\Rad(I+\p_1)=m$. Therefore,  the {\it Grothendieck's Non-vanishing Theorem} and \cite[Exercise 2.1.9]{BS} yield that $$H^2_{I}(R/\p_1)\cong H^2_{I+\p_1}(R/\p_1)=H^2_{\m}(R/\p_1)\neq 0.$$ Since, $\cd(I,R)=2$, the exact sequence $$0 \longrightarrow \p_1/\q_2 \longrightarrow R/\q_2 \longrightarrow R/\p_1 \longrightarrow 0$$ induces an exact sequence $$H^2_I(R/\q_2)\longrightarrow H^2_I(R/\p_1) \longrightarrow 0,$$ which implies that $H^2_I(R/\q_2)\neq 0$. So, $$\cd((I+\q_2)/\q_2, R/\q_2)=\cd(I,R/\q_2)\geq2.$$ Moreover, we have $(I+\q_2)/\q_2\subseteq \q_1/\q_2$ and hence $$\height (I+\q_2)/\q_2 \leq \height \q_1/\q_2=\height \q_1 - \height \q_2=2-1=1.$$
But, as $I\in\mathscr{I}(R)$ it follows from Lemma 3.1 that $$2\leq\cd(I,R/\q_2)=\height (I+\q_2)/\q_2\leq 1,$$ which is a contradiction.\qed\\

The following auxiliary lemmas are needed in the proof of Theorem 3.7.

 \begin{lem}
 \label{3.5}
 ${\rm (See}$ \cite[Proposition 4.7]{B2}$)$
The following statements are equivalent:
\begin{enumerate}[\upshape (i)]
  \item  $\mAss_R R= \mathfrak{A}(I,R) \cup\mathfrak{B}(I,R) \cup\mathfrak{D}(I,R)$, for any Noetherian complete local ring $(R,\m)$ and each $I\in\mathscr{I}(R)$,
  \item  $I\not\in\mathscr{I}(R)$, for any Noetherian complete local domain $(R,\m)$ of dimension $4$  and any ideal $I$ of $R$ with the property $\height I=2$.
  \end{enumerate}
\end{lem}
\qed\\

\begin{lem}
\label{3.6}
${\rm (See}$ \cite[Theorem 3.8]{B2}$)$
Let $I$ be an ideal of $R$ such that $$\mAss_R R=\mathfrak{A}(I,R)\cup \mathfrak{B}(I,R)\cup\mathfrak{D}(I,R).$$ Then $I\in\mathscr{I}(R)$.
 \end{lem}
\qed\\

The following theorem is  the second main result of this section.

\begin{thm}
\label{3.7}
Let $(R,\m)$ be a Noetherian complete local ring and $I$ be an ideal of $R$. Then the following statements are equivalent:
\begin{enumerate}[\upshape (i)]
  \item  $I\in\mathscr{I}(R)$.
  \item  $\mAss_R R= \mathfrak{A}(I,R) \cup\mathfrak{B}(I,R) \cup\mathfrak{D}(I,R)$.
  \end{enumerate}
 \end{thm}
\proof (i)$\Rightarrow$(ii) The assertion follows from Proposition 3.4 and Lemma 3.5. \\
(ii)$\Rightarrow$(i) The assertion holds by Lemma 3.6.\qed\\

In the final result of this paper, we apply Theorem 3.7 to the class of regular complete local rings.

\begin{cor}
\label{3.8}
Let $(S,\n)$ be a regular complete local ring. Then
  $$\mathscr{I}(S)=\{I\leq S\,\,:\,\,\Rad(I)=xS,\,\,{\rm for}\,\,{\rm some}\,\,x\in S\}\bigcup\{I\leq S\,\,:\,\, 0\leq\dim S/I\leq 1\}.$$
\end{cor}
\proof Set $$T:=\{I\leq S\,\,:\,\,\Rad(I)=xS,\,\,{\rm for}\,\,{\rm some}\,\,x\in S\}\bigcup\{I\leq S\,\,:\,\, 0\leq\dim S/I\leq 1\}$$ and assume that $I\in T$. Then, using \cite[Theorem 3.3.1]{BS} it is easy to see that $$\mAss_S S=\{0\}= \mathfrak{A}(I,S) \cup\mathfrak{B}(I,S) \cup\mathfrak{D}(I,S)$$and so it follows from Theorem 3.7 that $I\in \mathscr{I}(S)$. \\

Now, let $I\in \mathscr{I}(S)$. Then, by Theorem 3.7 we have $$\mAss_S S=\{0\}= \mathfrak{A}(I,S) \cup\mathfrak{B}(I,S) \cup\mathfrak{D}(I,S).$$
We consider the following cases:\\

\item Case 1. If $0\in \mathfrak{A}(I,S)$ then $I\subseteq 0$ or $I+0=S
$. So, $\Rad(I)=0S$ or $\Rad(I)=1_{_S}S$ and hence $I\in T$.\\

\item Case 2. If $0\in \mathfrak{B}(I,S)$ then we have $\cd(I,S)=1$ and so it follows from
\cite[Lemma 6.3.1 and Corollary 6.3.6]{BS} that  $\height \p=1$, for each $\p\in \mAss_S S/I$. But, as $$\Rad(I)=\bigcap_{\p\in \mAss_S S/I}\p$$ and $S$ is a unique factorization domain, it follows from \cite[Exercise 20.3]{Mat} that $\Rad(I)=xS$ for some $x\in S$ and hence $I\in T$.\\

\item Case 3. If $0\in \mathfrak{D}(I,S)$ then we have $0\leq\dim S/I\leq 1$ and so $I\in T$.\\

Now, we are ready to deduce that

$$\mathscr{I}(S)=T=\{I\leq S\,\,:\,\,\Rad(I)=xS,\,\,{\rm for}\,\,{\rm some}\,\,x\in S\}\bigcup\{I\leq S\,\,:\,\, 0\leq\dim S/I\leq 1\}.$$\qed\\

\subsection*{Acknowledgements}
The author would like to acknowledge his deep gratitude from
the referee for a very careful reading of the manuscript
and many valuable suggestions. He also, would like
to thank to School of Mathematics, Institute for Research in Fundamental
Sciences (IPM) for its financial support.


\end{document}